\documentclass[12pt]{article}
\usepackage{amsfonts,amssymb,version}

\def\F{{\cal F}} 
\def\mm{{\mathfrak m}}
\def\OO{{\mathcal O}}
\def\Q{\mathbb Q}
\def\Z{\mathbb Z}

\def\gr{\mathop{\rm Gr}\nolimits} 
\def\hn{\mathop{\rm HN}\nolimits}
\def\hnt{\mathop{\rm{HNT}}\nolimits}
\def\HNP{\mathop{\rm HNP}\nolimits}
\def\HOM{\mathop{\rm Hom}\nolimits}

\def\ker{\mathop{\rm ker}\nolimits}
\def\spec{\mathop{\rm Spec}\nolimits}

\def\toto{\stackrel{\to}{\scriptstyle \to}}

\let\hra\hookrightarrow
\let\ov\overline

\newtheorem{theorem}{Theorem}
\newtheorem{proposition}[theorem]{Proposition}
\newtheorem{lemma}[theorem]{Lemma}

\newtheorem{corollary}[theorem]{Corollary}

\def\rem{\refstepcounter{theorem}\paragraph{Remark \thetheorem}}

\begin{document}

\baselineskip=16pt

\centerline{\Large\bf Schematic Harder-Narasimhan Stratification}

\bigskip

\centerline{\bf Nitin Nitsure}

\bigskip

\begin{abstract} 
For any flat family of pure-dimensional coherent sheaves on a family
of projective schemes, the Harder-Narasimhan type (in the sense of 
Gieseker semistability) of its restriction to each fiber is known to 
vary semicontinuously on the parameter scheme of the family. This 
defines a stratification of the parameter scheme by locally closed
subsets, known as the Harder-Narasimhan stratification.

In this note, we show how to endow each Harder-Narasimhan stratum 
with the structure of a locally closed subscheme of the parameter
scheme, which enjoys the universal property that under any base change 
the pullback family admits a relative Harder-Narasimhan filtration
with a given Harder-Narasimhan type if and only if the base change 
factors through the  schematic stratum corresponding to that 
Harder-Narasimhan type.

The above schematic stratification induces a stacky stratification on 
the algebraic stack of pure-dimensional coherent sheaves. 
We deduce that coherent sheaves of a fixed Harder-Narasimhan type 
form an algebraic stack in the sense of Artin.
\end{abstract}

\bigskip

\centerline{2010 Math. Subj. Class. : 
14D20, 14D23, 14F05.}

\bigskip

\noindent{\large\bf 1. Introduction}

\medskip

Let $X$ be a projective scheme over a locally noetherian base scheme
$S$, with 
a chosen relatively ample line bundle $\OO_{X}(1)$.
Let $E$ be a coherent sheaf on $X$ which is flat over $S$, such that 
the restriction $E_s = E|_{X_s}$ of $E$ to the 
schematic fiber $X_s$ of $X$ over each $s \in S$ is a 
pure-dimensional sheaf of a fixed dimension $d\ge 0$.
For any $s\in S$, let $\hn(E_s)$ denote the Harder-Narasimhan type
of $E_s$ in the sense of Gieseker semistability. 
With respect to a certain natural partial order on the 
set $\hnt$ of all possible Harder-Narasimhan types $\tau$, the
Harder-Narasimhan function $s\mapsto \hn(E_s)$ is known to 
be upper semicontinuous on $S$. 

In this note, we prove that each level set $S^{\tau}(E)$ of the 
Harder-Narasimhan function has a natural structure of a locally
closed subscheme of $S$, with the following universal property: 
a morphism $T\to S$ factors via $S^{\tau}(E)$ if and only if
the pullback $E_t$ on $X_t$ for each $t\in T$ is of type $\tau$ and 
the pullback family $E_T$ on $X\times_ST$ admits a 
relative Harder-Narasimhan filtration, that is, a filtration 
$0\subset F_1\subset \ldots \subset F_{\ell} = E_T$ by coherent subsheaves 
such that the graded pieces $F_i/F_{i-1}$ are flat
over $T$, which for each $t\in T$ restricts to the 
Harder-Narasimhan filtration of $E_t$ on $X_t$. 

As a corollary, we deduce that sheaves of a fixed Harder-Narasimhan type
form an algebraic stack in the sense of Artin.

In Section 2 we recall the basic definitions and results 
of Harder-Narasimhan-Shatz that we need. In Section 3 we 
prove our main result 
(Theorem \ref{schematic Harder-Narasimhan stratification}), which
gives natural schematic structures  
on the Harder-Narasimhan strata. In Section 4, we show 
(Theorem \ref{algebraic stack}) that 
sheaves of a given Harder-Narasimhan type form an Artin algebraic stack.

This work had its origin in 
questions arising from the proposal of Leticia Brambila-Paz to
construct a moduli scheme for indecomposable unstable rank 2 vector 
bundles on a curve,  
fixing their Harder-Narasimhan type and the dimension of their 
vector space global endomorphisms. 
A construction of such a moduli scheme is given in [B-M-Ni], which uses
special cases of the results proved here.

\bigskip

\bigskip

\noindent{\large\bf 2. The Harder-Narasimhan filtration and stratification}

\medskip

Let $\Q[\lambda]$ be the polynomial ring in the variable $\lambda$. 
An element $f\in \Q[\lambda]$ is called a {\bf numerical polynomial} if 
$f(\Z) \subset \Z$. If a nonzero numerical polynomial  
$f$ has degree $d$, it can be uniquely expanded as 
$f = (r(f)/d!)\lambda^d + \mbox{ lower degree terms}$, where
$r(f) \in \Z$. If $f=0$ we put $r(f) =0$.
There is a {\bf total order} $\le$ on $\Q[\lambda]$ under which $f\le g$ 
if $f(m) \le g(m)$ for all sufficiently large integers $m$. 
Let the {\bf set of all Harder-Narasimhan types}, denoted by $\hnt$,
be the set consisting of all 
finite sequences $(f_1,\ldots, f_p)$ of numerical polynomials in 
$\Q[\lambda]$, where $p$ is allowed to vary
over all integers $\ge 1$,  
such that the following three conditions are satisfied.

\noindent{(1)} We have $0< f_1 < \ldots < f_p$ in $\Q[\lambda]$,\\
\noindent{(2)} the polynomials $f_i$ are all of the same degree,
say $d$, and \\
\noindent{(3)} the following inequalities are satisfied  
$$\frac{f_1}{r(f_1)}  > \frac{f_2 - f_1}{r(f_2) - r(f_1)} > \ldots >
\frac{f_p - f_{p-1}}{r(f_p) - r(f_{p-1})}.$$

Given any $x = (a,f)$ and $y= (b,g)$ in $\Z \times \Q[\lambda]$,
the {\bf segment joining $x$ and $y$} is the subset 
$\ov{xy} \subset \Z \times \Q[\lambda]$, 
consisting of all $(c,h)$ such that $(c,h)
= t(a,f) + (1-t)(b,g)$ for some $t\in \Q$ with $0\le t\le 1$.
For any 
$(f_1,\ldots, f_p)$ in $\hnt$, we define the corresponding
{\bf Harder-Narasimhan polygon} to be the subset 
$$\HNP( f_1,\ldots, f_p) \subset \Z \times \Q[\lambda]$$
which is the union of the segments 
$\ov{x_0x_1} \cup \ov{x_1x_2} \cup \ldots
\cup \ov{x_{p-1}x_p}$
where $x_0 = (0,0)$ and $x_i = (r(f_i), f_i)$ for $1\le i\le p$.

A point $(a,f)\in \Z \times \Q[\lambda]$ is said to {\bf lie under}
another point $(b,g )\in \Z \times \Q[\lambda]$ if 
$a=b$ in $\Z$ and $f\le g$ in $\Q[\lambda]$. 
A point $(a,f)\in \Z \times \Q[\lambda]$ is said to {\bf lie under the 
polygon} $\HNP(g_1,\ldots, g_q)$ if there exists some $(b,g) \in 
\HNP( g_1,\ldots, g_q)$ such that the point 
$(a,f)$ lies under the point $(b,g)$.
There is a {\bf partial order} $\le$ on $\hnt$, under which
$(f_1,\ldots, f_p) \le (g_1,\ldots, g_q)$
if for each $1\le i \le p$, the point 
$(r(f_i), f_i)$
lies under $\HNP(g_1,\ldots, g_q)$.

With the above numerical preliminaries, we now briefly recall
the theory of Harder-Narasimhan-Shatz for filtrations and stratifications.
For an exposition 
the reader can see, for example, the book of Huybrechts and Lehn [Hu-Le]. 

Let $Y$ be a projective scheme over a field $k$, together
with an ample line bundle $\OO_Y(1)$. For any coherent sheaf 
$E$ on $Y$, we denote by $P(E)\in \Q[\lambda]$ the resulting 
Hilbert polynomial of $E$, defined by 
$P(E)(m) = \sum_i (-1)^i\dim_k H^i(Y, E(m))$.
This is a numerical polynomial, and the integer
$r(P(E))$ is denoted simply by $r(E)$ (this is called the {\bf rank}
of $E$). For a nonzero $E$, the degree $d$ of $P(E)$ equals the dimension of 
the support of $E$. A coherent sheaf $E$ on $Y$ is said to be 
{\bf pure-dimensional} of dimension $d$ 
if for every open subscheme $U \subset Y$ and 
nonzero coherent subsheaf $F\subset E|_U$, the support of $F$ is of the 
same dimension $d$. A pure-dimensional coherent sheaf $E$ is called 
{\bf semistable} (in the sense of Gieseker)
w.r.t. $\OO_Y(1)$ if for all coherent subsheaves $F\subset E$, we have
$r(E)P(F) \le r(F)P(E)$. 

When $E$ is of pure dimension $d\ge 0$ but not necessarily semistable, 
it admits a unique strictly increasing filtration 
$0= \hn_0(E) \subset\hn_1(E)\subset \ldots \subset \hn_{\ell}(E) = E$
by coherent subsheaves $\hn_i(E)$ such that for each $1\le i\le \ell$, 
the graded piece $\gr_i(E) = \hn_i(E)/\hn_{i-1}(E)$ is semistable of
pure dimension $d$, and the inequalities
$$\frac{P(\gr_1(E))}{r(\gr_1(E))} > \ldots > 
\frac{P(\gr_{\ell}(E))}{r(\gr_{\ell}(E))}$$
hold. This filtration is 
called the {\bf Harder-Narasimhan filtration} of $E$
(in the sense of Gieseker semistability). 
The first step $\hn_1(E)$ is called the {\bf maximal destabilizing
subsheaf} of $E$.
The integer $\ell$ (also written as $\ell(E)$) is called as the 
{\bf length} of the Harder-Narasimhan filtration of $E$. 
In these terms, a nonzero pure-dimensional sheaf is semistable 
if and only if its Harder-Narasimhan filtration is of length $\ell(E) = 1$.
The ordered $\ell(E)$-tuple 
$$\hn(E) = (P(\hn_1(E)), \ldots ,  P(\hn_{\ell}(E)))\in \hnt$$
is called the {\bf Harder-Narasimhan type} of $E$.

In his path-breaking paper [Sh], S.S. Shatz addressed the 
question of the variation of the Harder-Narasimhan type in a family. 
The set-up for this is as follows. 
Let $S$ be a locally noetherian scheme, and let 
$\pi: X \to S$ be a projective scheme over $S$, 
with a relatively ample line bundle $\OO_{X}(1)$. 
Let $E$ be a coherent sheaf of $\OO_{X}$-modules 
which is flat over $S$ such each 
restriction $E_s$ 
to the schematic fiber $X_s = \pi^{-1}(s)$ is pure-dimensional.
The {\bf Harder-Narasimhan function of $E$} is the function
$$|S| \to \hnt : s\mapsto \hn(E_s)$$
where $|S|$ denotes the underlying topological space  
of the scheme $S$. Shatz proved in [Sh] that $\hn(E_s)$ is 
upper-semicontinuous w.r.t. the partial order $\le$ on
$\hnt$ defined above (actually, HN-filtrations in the sense of 
$\mu$-semistability rather than Gieseker semistability are
considered in [Sh], but the proofs in the Gieseker semistability 
case are similar with obvious changes).

\rem\label{strata are locally closed} In particular, 
for any $\tau \in \hnt$, the corresponding level set 
$$|S|^{\tau}(E) =\{ s \in |S| \mbox{ such that } \hn(E_s) = \tau \}$$ 
is locally closed in $|S|$, 
the subset $|S|^{\le \tau}(E)  =\bigcup_{\alpha \le \tau}
|S|^{\alpha}(E) \subset |S|$ is
open in $|S|$, and $|S|^{\tau}(E)$ is closed in $|S|^{\le \tau}(E)$.

\rem\label{why the quotient is again pure} 
If $(f_1,\ldots, f_p)\in \hnt$, then $(f_2-f_1, \ldots, f_p-f_1)$
is again in $\hnt$. Let $E$ be pure-dimensional on $Y$ with 
$\hn(E) \le (f_1,\ldots, f_p)\in \hnt$.
If $E'\subset E$ is a coherent subsheaf with $P(E') = f_1$,
then we must have $\hn_1(E) = E'$, that is, such an $E'$
is automatically the maximal destabilizing subsheaf of $E$.
The quotient $E'' = E/E'$ is pure-dimensional, with 
$\hn(E'') \le (f_2-f_1, \ldots, f_p-f_1)$. Moreover, we have
$\HOM_Y(E', E'') = 0$.

\rem\label{HN base changes}
If $(Y,\OO_Y(1))$ is a projective scheme over a field $k$ and if
$K$ any extension field of $k$, then a coherent sheaf $E$ on $Y$ is
semistable  w.r.t. $\OO_Y(1)$ if and only if its base-change 
$E_K = E\otimes_kK$ to $Y_K$ is semistable  w.r.t. $\OO_{Y_K}(1)
= \OO_Y(1)\otimes_kK$. Consequently, if $E$ is any
pure-dimensional sheaf on $Y$ then the Harder-Narasimhan filtration
$\hn_i( E_K)$ is just the pullback $\hn_i( E)\otimes_kK$
of the Harder-Narasimhan filtration of $E$.

\bigskip

\bigskip

\noindent{\large\bf 3. Scheme structures on HN strata}

\bigskip

For basic facts that we need from 
Grothendieck's theory of Quot schemes and their deformation theory,
the reader can consult [Hu-Le], [F-G], [Ni 1] and [Ni 2].

\begin{lemma}\label{closed embedding}
A morphism $f: T\to S$ between locally noetherian schemes
is a closed embedding if (and only if) 
$f$ is proper, injective, unramified and induces an 
isomorphism $k(f(t)) \to k(t)$ of residue fields for all $t\in T$.
\end{lemma}

\medskip \noindent{{\bf Proof}} 
Note that $f(T)$ is closed in $S$, and $f_*\OO_T$ is coherent. 
It only remains to show that the homomorphism 
$f^{\#} : \OO_S\to f_*\OO_T$ is surjective. 
It is enough to show it stalk-wise at all points
of $f(T)$, so we can assume that 
$S = \spec A$ where $A$ is a noetherian local ring. Then by 
finiteness and injectivity of $f$, we have $T = \spec B$ where 
$B$ is a finite local $A$ algebra, and 
$f^{\#} : A\to B$ is a local homomorphism which by assumption
induces an isomorphism $A/\mm_A \to B/\mm_B$ on the residue fields, and 
$\mm_B = \mm_A B$ by assumption of unramifiedness. Hence it follows by the 
Nakayama lemma that $f^{\#} : A\to B$ is surjective.
$\Box$

\bigskip

Let $X$ be a projective scheme over a locally noetherian base scheme
$S$, with 
a chosen relatively ample line bundle $\OO_{X}(1)$.
Let $E$ be a coherent sheaf on $X$ which is flat over $S$, such that 
the restriction $E_s = E|_{X_s}$ of $E$ to the 
schematic fiber $X_s$ of $X$ over each $s \in S$ is a nonzero 
pure-dimensional sheaf of a fixed HN type $\tau = (f_1,\ldots, f_{\ell})$. 
A {\bf relative Harder-Narasimhan filtration} of $E$ 
is a filtration $0= E_0 \subset E_1\subset \ldots \subset E_{\ell} = E$
by coherent subsheaves on $X$, such that for each $i$ 
with $1\le i \le \ell$, the quotient $E_i/E_{i-1}$ is flat $S$,
and for each $s\in S$ this filtration restricts to give the 
Harder-Narasimhan filtration $\hn_i(E_s)$ of $E_s$. 

We now come to the main result of this note.

\begin{theorem}\label{schematic Harder-Narasimhan stratification}
{\bf (Main Theorem)}
Let $X$ be a projective scheme over a locally noetherian scheme $S$, with 
a relatively ample line bundle $\OO_{X}(1)$.
Let $E$ be a coherent sheaf on $X$ which is flat over $S$, such that 
the restriction $E_s$ is a pure-dimensional sheaf on $X_s$ 
for each $s \in S$. Let $\tau = (f_1,\ldots, f_{\ell})\in \hnt$.
Then we have the following.

{\noindent (1)} Each Harder-Narasimhan stratum 
$|S|^{\tau}(E)$ of $E$ has a unique structure of a locally
closed subscheme $S^{\tau}(E)$ of $S$, with the following universal
property: a morphism $T\to S$ factors via $S^{\tau}(E)$ if and only if
the pullback $E_T$ on $X\times_ST$ 
admits a relative Harder-Narasimhan filtration of type $\tau$.

{\noindent (2)} A relative Harder-Narasimhan 
filtration on $E$, if it exists, is unique.

{\noindent (3)} For any morphism $f: T\to S$ of locally noetherian schemes,
the schematic stratum  $T^{\tau}(E_T)\subset T$ for $E_T$ equals the 
schematic inverse image of $S^{\tau}(E)$ under $f$. 
\end{theorem}

\medskip \noindent{{\bf Proof}} 
If $\ell =1$, then we take $S^{\tau}(E)$ to be the open
subscheme of $S$ consisting of all $s$ such that $E_s$
is semistable with Hilbert polynomial $f_1$. We now argue
by induction on $\ell \ge 2$. 
By Remark \ref{strata are locally closed}, 
all $s$ with $\hn(E_s) \le \tau$ form an open subset 
$|S|^{\le \tau}(E)$ of $S$, and $|S|^{\tau}(E)$ is a closed
subset of $|S|^{\le \tau}(E)$. We give 
$|S|^{\le \tau}(E)$ the unique structure of an open subscheme of $S$,
which we denote by $S^{\le \tau}(E)$. 
In what follows we 
will give the closed subset $|S|^{\tau}(E)$
a particular structure of a closed subscheme of $S^{\le \tau}(E)$, 
which has the desired universal property.

Let $X^{\le \tau}$ be the inverse image of 
$S^{\le \tau} = S^{\le \tau}(E)$ in $X$,
and let $\OO_{X^{\le \tau}}(1)$ and $E^{\le \tau}$ be the restrictions  
of $\OO_X(1)$ and $E$ to $X^{\le \tau}$.
Consider the relative Quot scheme 
$$Q = Quot_{E^{\le \tau}/X^{\le \tau}/S^{\le \tau}}^{f_{\ell} - f_1 ,\,\OO_{X^{\le \tau}}(1)}$$ 
with projection $\pi : Q \to S^{\le \tau}$. 
Then $\pi$ is projective, hence proper.

Let $q\in Q$ represent a quotient $q' : E_q \to \F$ on $X_q$. 
Then $\ker(q') = \hn_1(E_q)$ by Remark \ref{why the quotient is again pure}.
If $q\mapsto s\in S^{\le \tau}$, then by Remark \ref{HN base changes}
the quotient $q'$ is the pullback of the quotient 
$E_s\to E_s/\hn_1(E_s)$ which is defined over $X_s$. 
Hence the residue field extension $k(s) \to k(q)$ is trivial. 
By the uniqueness of $\hn_1(E_s)$, there exists 
at most one such $q$ over $s$.
The fiber of $\pi: Q \to S^{\le \tau}$ over $s$ is the Quot scheme 
$$\pi^{-1}(s) = Quot_{E_s/X_s/k(s)}^{f_{\ell} - f_1 ,\, \OO_{X_s}(1)}.$$
By a standard fact in the deformation theory for Quot schemes
(see, for example, Theorem 3.11.(2) in [Ni 2]), 
its tangent space at $q$ is given by
$$T_q(\pi^{-1}(s))= \HOM_{X_q}(\ker(q'),\,E_q/\ker(q'))
= \HOM_{X_q}(\hn_1(E_q), E_q/\hn_1(E_q))$$ 
which is zero by 
Remark \ref{why the quotient is again pure}.
Hence $\pi: Q \to S^{\le \tau}$ is unramified. 

It now follows by Lemma \ref{closed embedding} 
that $\pi: Q\to S^{\le \tau}$ is a closed imbedding. 

Now consider the universal quotient sheaf 
$E_Q \to E''$ on $X_Q = X\times_SQ$. 
By Remark \ref{why the quotient is again pure},
for all $q \in Q$ the sheaf $E''_q$ on $X_q$ is
pure-dimensional, with 
$$\hn(E''_q) \le \tau'' = (f_2-f_1,\ldots, f_{\ell} - f_1).$$
In particular, we have $Q^{\le \tau''}(E'') = Q$.
The Harder-Narasimhan type $\tau''$ has length $\ell -1$, 
hence by induction on the length, the closed subset 
$|Q|^{\tau''}(E'')$ of $Q$ has the structure of a closed 
subscheme $Q^{\tau''}(E'') \subset Q$ which has the desired
universal property for $E''$. We regard $Q$ as a closed subscheme of 
$S^{\le \tau}$ via $\pi$, and we finally 
define the closed subscheme $S^{\tau}(E)\subset S^{\le \tau}$ by
putting  
$$S^{\tau}(E) = Q^{\tau''}(E'') \subset Q \subset S^{\le \tau}.$$

We now show that $S^{\tau}(E)$ so defined has the desired universal 
property. As $\tau''$ has length $\ell -1$, 
by induction on the length, the restriction $E''_{S^{\tau}(E)}$ 
of $E''$ to $X\times_S S^{\tau}(E)$ 
has a unique relative Harder-Narasimhan filtration
$$0 \subset E''_1 \subset \ldots  \subset E''_{\ell -1} = E''_{S^{\tau}(E)}$$
with Harder-Narasimhan type $\tau''$. 
For $2\le i \le \ell$, let $E_i$ be the inverse image of $E''_{i-1}$ 
under the restriction of universal quotient
$E_Q \to E''$ to $X_{S^{\tau}(E)}$. 
This defines a relative Harder-Narasimhan filtration 
$0 \subset E_1 \subset \ldots  \subset E_{\ell} = E_{S^{\tau}(E)}$ 
of $E_{S^{\tau}(E)}$ over the base $S^{\tau}(E)$. In particular, 
if a morphism $T\to S$ factors via $S^{\tau}(E)$ then the  
pullback of this filtration gives a 
relative Harder-Narasimhan filtration over $T$. 

Conversely, let $f: T\to S$ be a morphism such that the pullback
$E_T$ on $X_T$ has a relative Harder-Narasimhan filtration 
$0= F_0 \subset F_1 \subset \ldots \subset F_{\ell} =E_T$
of type $\tau$. 
The quotient $E_T \to E_T/F_1$ has Hilbert polynomial $f_{\ell} - f_1$ over
all $t\in T$, so by the universal property of the Quot scheme $Q$, 
the morphism $T\to S$ factors via $Q\hra S$,
inducing a morphism $f' : T\to Q$.  
By Remark \ref{why the quotient is again pure}, the restriction of 
$E_T/F_1$ is pure-dimensional on $X_t$ for all $t\in T$, and 
$0 = (F_1/F_1) \subset (F_2/F_1) \subset \ldots \subset
(F_{\ell}/F_1) = E_T/F_1$
is a relative Harder-Narasimhan filtration of 
$E_T/F_1 = (f')^*(E'')$ over the base $T$, 
with type $\tau''$ which has length $\ell -1$. 
Hence by induction,   
$f': T\to Q$ factors via $Q^{\tau''}(E'') = S^{\tau}(E)$, as
desired. This completes the proof of (1).

Next we show the uniqueness of 
a relative Harder-Narasimhan filtration
$0= E_0 \subset E_1\subset \ldots \subset E_{\ell} = E$ 
over a base $S$, assuming such a filtration exists. As 
$|S|^{\tau}(E) = |S|$, we at least have $S = S^{\le \tau}(E)$. 
With notation as above, we have shown that
$\pi: Q\to S^{\le \tau}$ is a closed imbedding, therefore 
$\pi$ admits at most one global section, which shows that $E_1$ is unique. 
By inductive assumption on $\ell$, the quotient family
$E/E_1$ admits a unique relative Harder-Narasimhan filtration
$F_j$, so defining $E_i\subset E$ to be the inverse image of  
$F_{i-1}$ under $E\to E/E_1$ for $2\le i\le \ell$, we see that 
$0= E_0 \subset E_1\subset \ldots \subset E_{\ell} = E$ is the only 
possible relative Harder-Narasimhan filtration on $E$. 
This proves the statement (2).

The base-change property (3) for the schematic strata 
is a direct consequence of the universal property (1). 
This completes the proof of the theorem.
\hfill$\Box$

\bigskip

The following immediate implication of Theorem 
\ref{schematic Harder-Narasimhan stratification}
shows that when the HN type is constant over a reduced base scheme, 
the outcome is as nice as can be expected.

\begin{corollary}{\bf (Case of constant HN type over a reduced base)}
Let $X$ be a projective scheme over a locally noetherian base scheme
$S$, with a chosen relatively ample line bundle $\OO_{X}(1)$.
Let $E$ be a coherent sheaf on $X$ which is flat over $S$, such that 
the restriction $E_s = E|_{X_s}$ of $E$ to the 
schematic fiber $X_s$ of $X$ over each $s \in S$ is a 
pure-dimensional sheaf of a fixed Harder-Narasimhan type 
$\tau \in \hnt$. Suppose moreover that $S$ is 
reduced. Then $S = S^{\tau}$, that is, 
$E$ admits a unique relative Harder-Narasimhan filtration.
\end{corollary}

\bigskip

\noindent{\large\bf 4. Moduli stack $Coh_{X/S}^{\tau}$}

\medskip

For basic terminology and conventions about stacks, we will follow
the book [L-MB] by Laumon and Moret-Bailly.
In what follows, $X$ will 
be a projective scheme over a locally noetherian base 
scheme $S$, with 
a chosen relatively ample line bundle $\OO_{X}(1)$, and
$\tau = (f_1,\ldots,f_{\ell})$ will be any element of  $\hnt$.

Let $Coh_{X/S}$ denote the Artin algebraic stack over $S$ 
of all flat families of coherent sheaves on $X/S$ (see [L-MB] 2.4.4).
In any such family, pure-dimensionality of all 
restriction to fibers is an open condition on the parameter scheme, 
pure-dimensionality is preserved by arbitrary base changes, 
and the base-change under a surjection is pure-dimensional 
on all fibers if and only if the original is so. Hence
pure-dimensional coherent sheaves form an open algebraic substack
$Coh_{X/S}^{pure} \subset Coh_{X/S}$. 

We will define the moduli stack $Coh_{X/S}^{\tau}$ 
of pure-dimensional coherent sheaves of type $\tau$ as
a strictly full sub $S$-groupoid of $Coh_{X/S}^{pure}$,
as follows. For any $S$-scheme $T$, 
we say that an object $E\in Coh_{X/S}^{pure}(T)$
lies in $Coh_{X/S}^{\tau}(T)$ if and only if 
$E$ admits a relative Harder-Narasimhan filtration 
with constant type $\tau$. This is clearly closed under
pullbacks $f^* : Coh_{X/S}(T)\to Coh_{X/S}(T')$ 
for all $S$-morphisms $f: T'\to T$. 

To prove that the $S$-groupoid $Coh_{X/S}^{\tau}$ 
thus defined is a stack,
we need the following property of effective descent.

\begin{lemma}\label{fpqc descent of Harder-Narasimhan} 
Let $T$ be an $S$-scheme and let $E$ be an object of $Coh_{X/S}(T)$. 
Let $f: T'\to T$ be a faithfully flat quasi-compact morphism. 
If the pullback $f^*E$ is in $Coh_{X/S}^{\tau}(T')$,
then $E$ is in $Coh_{X/S}^{\tau}(T)$. 
\end{lemma}

\medskip \noindent{{\bf Proof}} 
Each $E_t$, where $t\in T$,  
is pure-dimensional with Harder-Narasimhan type $\tau$, 
as its pullback $E_{t'}$ is so for any $t'\in T'$ with $t'\mapsto t$,
and as $T'\to T$ is surjective.
It now only remains to construct a relative 
Harder-Narasimhan filtration of $E$. This we do by showing that 
the relative
Harder-Narasimhan filtration $(F_i)$ of the pullback $E_{T'}$
descends under $T'\to T$.

Let $T'' = T'\times_TT'$ with projections $\pi_1,\pi_2: T'' \toto T'$. 
By Grothendieck's result on effective fpqc descent for quasicoherent
subsheaves of the pullback of a quasicoherent sheaf, 
to show that the filtration descends to $T$ 
we just have to show that the pullbacks of the filtration
under the two projections $\pi_1,\pi_2: T'' \toto T'$ are identical.
But note that we have an identification 
$\pi_1^*(E_{T'}) = \pi_2^*(E_{T'}) = E_{T''}$,
under which the pullbacks $\pi_1^*(F_i)$ and $\pi_2^*(F_i)$
are relative Harder-Narasimhan filtrations of $E_{T''}$. 
Hence these filtrations coincide by 
Theorem \ref{schematic Harder-Narasimhan stratification}.
\hfill$\Box$

\begin{theorem}\label{algebraic stack}
Let $X$ be a projective scheme over a locally noetherian scheme $S$, with 
a relatively ample line bundle $\OO_{X}(1)$.
Let $\tau$ be any Harder-Narasimhan type. Then all flat families of
pure-dimensional coherent sheaves on $X/S$ with fixed 
Harder-Narasimhan type $\tau$ form an 
algebraic stack $Coh_{X/S}^{\tau}$ over $S$, which is a 
locally closed substack of the algebraic stack $Coh_{X/S}$ of 
all flat families of coherent sheaves on $X/S$. 
\end{theorem}

\noindent{{\bf Proof}} The inclusion $1$-morphism of $S$-groupoids 
$\theta : Coh_{X/S}^{\tau} \hra Coh_{X/S}^{pure}$
is fully faithful. Hence $Coh_{X/S}^{\tau}$ is a pre-stack over $S$.
By Lemma \ref{fpqc descent of Harder-Narasimhan}, the pre-stack
$Coh_{X/S}^{\tau}$ satisfies effective fpqc descent, so it is a stack
over $S$. We next prove that it is algebraic.

Given any $E$ in $Coh_{X/S}^{pure}(T)$, 
let $T^{\tau}(E)\subset T$ be the corresponding 
schematic Harder-Narasimhan stratum as given 
by Theorem \ref{schematic Harder-Narasimhan stratification}.
Let $[E] : T \to  Coh_{X/S}^{pure}$ be the classifying $1$-morphism of $E$. 
By Theorem \ref{schematic Harder-Narasimhan stratification}
we have a natural isomorphism 
$$T \,\times_{[E], Coh_{X/S}^{pure}, \theta}\,Coh_{X/S}^{\tau} 
\, \cong \, T^{\tau}(E)$$
of $S$-groupoids, under which the projection of 
the fibered product to $T$ corresponds to the imbedding of 
$T^{\tau}(E)$ as a locally closed subscheme in $T$. 

This shows the inclusion $1$-morphism 
$\theta : Coh_{X/S}^{\tau} \hra Coh_{X/S}$
of stacks is a representable locally closed imbedding. 
Hence $Coh_{X/S}^{\tau}$ is an algebraic stack over $S$,
which is a locally closed substack of $Coh_{X/S}$.
\hfill$\Box$

\bigskip

We now come to the question of quasi-projectivity of $Coh_{X/S}^{\tau}$.
For a given $X/S$, $\OO_{X}(1)$ and 
$\tau = (f_1,\ldots,f_{\ell})$,
consider the following {\bf boundedness condition} {\bf (*)}.

\medskip

\noindent{{\bf (*)}} {\it There exists a natural number
$N$ such that for any morphism $\spec K\to S$ where $K$ is a field 
and any semistable coherent sheaf $F$ on the base-change $X_K$  
whose Hilbert polynomial is equal to $f_i$ for any $1\le i\le \ell$,
the sheaf $F(N) = F\otimes \OO_{X_K}(N)$ is generated by global sections,
and all its cohomology groups $H^j(X_K, F(N))$ vanish for $j\ge 1$.
}

\medskip

By the boundedness theorems of Maruyama-Simpson [Si] and Langer [La], 
the condition {\bf (*)} is indeed satisfied in many cases of interest,
for example, when $S$ is of finite type over an algebraically
closed field $k$ of arbitrary characteristic.

\begin{proposition}{\bf (Quasi-projectivity of $Coh_{X/S}^{\tau}$)}
If the above boundedness condition {\bf (*)} is satisfied, 
then the stack $Coh_{X/S}^{\tau}$ admits an atlas 
$U\to Coh_{X/S}^{\tau}$ such that $U$ is a quasi-projective scheme over $S$. 
\end{proposition}

\medskip \noindent{{\bf Proof}}   
If a coherent sheaf $E$ is of type $\tau$ 
on $X_K$ for an $S$-field $K$, then by 
{\bf (*)}, $E$ is a quotient of $\OO_{X_K}(-N)^{f_{\ell}(N)}$. 
Let $Q$ be the relative Quot scheme over $S$ which 
parameterizes all coherent quotient sheaves of $\OO_{X}(-N)^{f_{\ell}(N)}$
on fibers of $X/S$, with fixed Hilbert polynomial $f_{\ell}$.
Let $E$ be the universal quotient sheaf on $X\times_SQ$.
Let $Q_o \subset Q$ be the open subscheme 
consisting of all $q\in Q$ satisfying the conditions that $E_q$  
is pure-dimensional, $E_q(N)$ is generated by global sections, 
the map $H^0(X_q, \OO_{X_q}^{f_{\ell}(N)}) \to H^0(X_q, E_q(N))$ 
induced by $q$ is an
isomorphism, and $H^i(X_q, E_q(N))=0$ for all $i\ge 1$ (each of
these conditions is an open condition). 

Let $E_o$ be the restriction of $E$ to $X\times_S Q_o$.
Let $Q_o^{\tau}(E_o)$ be the locally closed subscheme of $Q_o$ 
corresponding to the Harder-Narasimhan type $\tau$, given by 
Theorem \ref{schematic Harder-Narasimhan stratification}. 
The classifying $1$-morphism $[E_o]: Q_o^{\tau}(E_o)\to Coh_{X/S}^{\tau}$
of $E_o$ is an atlas for $Coh_{X/S}^{\tau}$
(that is, $[E_o]$ is a representable smooth surjection), 
as follows from the proof of 
Theorem 4.6.2.1 in Laumon and Moret-Bailly [L-MB]. As $Q$ is
projective over $S$, its locally 
closed subscheme $U$ is quasi-projective over $S$, as desired. 
\hfill$\Box$

\bigskip

\newpage

\centerline{\large \bf References}

\bigskip

\parindent 0pt

\parskip 4pt


[B-M-Ni] Brambila-Paz L., Mata, O. and Nitsure, N. :  
Endomorphisms and moduli for unstable bundles. (Preprint, 2009).

[F-G] Fantechi, B. and G\"ottsche, L. : Local properties and 
Hilbert schemes of points. Part 3 of 
{\it `Fundamental Algebraic Geometry -- Grothendieck's FGA
Explained.'},  Fantechi et al, Math. Surveys and Monographs Vol.
123, American Math. Soc. (2005).

[H-N] Harder, G. and Narasimhan, M.S. : 
On the cohomology groups of moduli spaces of vector bundles on curves.  
Math. Ann. 212 (1974/75), 215--248.

[Hu-Le] Huybrechts, D. and Lehn, M. : {\it `The Geometry of Moduli
Spaces of Sheaves.'} Aspects of Mathematics 31, Vieweg (1997).

[La] Langer, A. : Semistable sheaves in positive characteristic.  
Ann. of Math. (2) 159  (2004),  251--276. 

[L-MB] Laumon, G. and Moret-Bailly, L. : {\it `Champs alg\'ebriques.'}
Springer (2000).

[Ni 1] Nitsure, N. : Construction of Hilbert and Quot schemes. Part
2 of {\it `Fundamental Algebraic Geometry -- Grothendieck's FGA
Explained.'},  Fantechi et al, Math. Surveys and Monographs Vol.
123, American Math. Soc. (2005).

[Ni 2] Nitsure, N. : Deformation theory for vector bundles. Chapter 5
of {\it Moduli Spaces and Vector Bundles} 
(edited by Brambila-Paz, Bradlow, Garcia-Prada and Ramanan), 
London Math. Soc. Lect. Notes 359,  
Cambridge Univ. Press (2009).

[Sh] Shatz, S.S. : The decomposition and specialization of
algebraic families of vector bundles. Compositio Math. 35
(1977), no. 2, 163--187.

[Si] Simpson, C. Moduli of representations of the fundamental 
group of a smooth projective variety -I.  
~~Publ. Math. I.H.E.S. 79  (1994), 47--129.


\bigskip

\bigskip

\parskip 1pt

{\small 

School of Mathematics, \hfill {\tt nitsure@math.tifr.res.in}

Tata Institute of Fundamental Research, 


Homi Bhabha Road,

Mumbai 400 005,


INDIA.

\bigskip

\centerline{04 Sep 2009}

}

\end{document}